\documentclass{amsart}
\usepackage{amsmath,amsthm,amssymb,amsfonts}
\usepackage[colorlinks=true]{hyperref}

\hypersetup{urlcolor=blue, citecolor=blue}

\let\MR\mr

\def\doi#1{   {\href{http://dx.doi.org/#1}
   {{\mdseries\ttfamily DOI}}}}

    \newcommand{\be}{\beta}
    
  \newcommand{\ep}{\varepsilon}

    \newcommand{\Ga}{\Gamma}
\newcommand{\R}{\mathbb{R}}

\newcommand{\pt}{\partial_t}\newcommand{\pa}{\partial}
\newcommand{\les}{{\lesssim}}
\newcommand{\beeq}{\begin{equation}}\newcommand{\eneq}{\end{equation}}

\newenvironment{prf}{\noindent {\bf Proof.} }{\endprf\par}
\def \endprf{\hfill  {\vrule height6pt width6pt depth0pt}\medskip}

\numberwithin{equation}{section}

\newcommand{\gm}{\mathfrak{g}}

\def\<{\langle}             \def\>{\rangle}
\def\({\left(}                 \def\){\right)}

\newtheorem{theorem}{Theorem}[section]
\newtheorem{corollary}[theorem]{Corollary}

\newtheorem{lemma}[theorem]{Lemma}

\theoremstyle{definition}

\newtheorem{remark}[theorem]{Remark}

\theoremstyle{definition}
\title[Global existence for semilinear damped wave equations]
      {Global existence for the 3-D semilinear damped wave equations in the scattering case}

\author{Yige Bai}\address{School of Mathematical Sciences\\                Zhejiang University\\                Hangzhou 310027, P. R. China}\email{xqsyyg@163.com}

\author{Mengyun Liu}
\address{School of Mathematical Sciences\\                Zhejiang University\\                Hangzhou 310027, P. R. China}
\email{mengyunliu@zju.edu.cn}

\date{\today}
\dedicatory{} \commby{}

\begin{document}

\begin{abstract}
We study the global existence of solutions to semilinear damped wave equations in the scattering case with derivative power-type nonlinearity on (1+3) dimensional nontrapping asymptotically Euclidean manifolds. The main idea is to exploit local energy estimate, together with local existence to convert the parameter $\mu$ to small one.
\end{abstract}


\maketitle
%\tableofcontents
%%% Section 1 %%%
\section{Introduction}

In this work, we consider the global existence of solutions for the Cauchy problem of the damped semilinear wave equation
\begin{align}
\label{1.1}
&u_{tt} - \Delta_{\gm}u + \frac{\mu \pt{u}}{(1 + t)^\be}  = |\pt{u}|^p,
\end{align}
with initial data
\begin{align}
\label{ea2}
&u(0,x)=f(x),  \pt{u}(0,x)=g(x),
\end{align}
where $\be > 1$, and $\mu$, $\be \in \R$. Here we consider the nontrapping asymptotically Euclidean manifolds $(\R^{n}, \gm)$, and
$$\Delta_{\gm} = \sum_{i,j=1}^{n}\sqrt{|\gm|}^{-1}\partial_{i}g^{ij}\sqrt{|\gm|}\partial_{j},$$
\begin{align}
\label{ea6}
&\gm=g_{0}+g_1(r)+g_2(x), \gm ~is ~~nontrapping,
\end{align}
where $g^{ij}_{0}= \delta^{ij}$, $g_1$ and $g_2$ are of the form $g_{jk}dx^jdx^k$, $\langle x\rangle = \sqrt{1+x^{2}}$,
\begin{align}
\label{ea7}
&|\nabla_x^a g_{i,jk}|\les_a\langle x\rangle^{-|a|-\rho_i},i=1,2,\rho=\min(\rho_1,\rho_2-1),\rho_1<\rho_2,
\end{align}
where $\rho_{1}>0, \rho_{2}>1$ are fixed and $(g^{ij}(t,x))$ denotes the inverse matrix of $(g_{ij}(t,x))$. We assume the first perturbation $g_1$ is radial.\par
When $\gm = g_{0}$, and $\mu =0$, Glassey made a conjecture that the critical power $p$ for the problem to admit global solutions with small, smooth initial data with compact support is
$$p_c:=1 + \frac{2}{n-1}$$
in Glassey \cite{g} (see also Schaeffer \cite{scha} and Rammaha \cite{ram}), where n is the spatial dimension.  The conjecture was verified for $n=2,3$ for general data ( Hidano and Tsutaya \cite{hk} and Tzvetkov \cite{nt} independently, as well as the radical case in Sideris \cite{sideris} for $n=3$). Zhou \cite{zhou} obtained the blow up results for $n\geq 4$, when $p\leq p_c$ and Hidano, Wang, Yokoyama \cite{hyw} obtained the global existence for $p > p_{c}$ when $n \geq 4$ under the radically symmetric assumption. For the results on the Glassey conjecture on certain asymptotically flat manifolds ($g_{1}, g_{2} \neq 0$): $n=3$ was proved by Wang \cite{W15} on certain small space time perturbation of the flat metric, as well as the nontrapping asymptotically Euclidean manifolds. Moreover, in Wang \cite{W15} the Glassey conjecture was proved when $n \geq 3$ in radical case on radical asymptotically flat manifolds. See also Wang \cite {W15j} for related results on exterior domain with nontrapping obstacles.\par
When $\gm = g_{0}$, $\mu > 0$, for the corresponding linear problem
\beeq
\label{1.2}
\Box u + \mu\frac{u_{t}}{(1+t)^{\beta}} = 0
\eneq
we say that the damped term is ``scattering'' when $\beta \in (1, +\infty)$ since the solution behaves like that of wave equation (see, e.g., Wirth \cite{wirth1} \cite{wirth2} \cite{wirth3} for the classifications of (\ref{1.2})). \par
When $\gm = g_{0}$ and $\mu > 0$, Lai and Takamura \cite{ltm} proved that (\ref{1.1}) blows up at finite time when $1 <p \leq p_{c}$ and $n \geq 1$. In view of the results of \cite{W15} and the ``scattering" damping term, it is natural to expect (\ref{1.1}) would admit global solution for any $\mu \in \R$ when $n=3$ and $p >2$. In following Theorem, we prove that it is the case. See Theorem \ref{4.1} for a more precise statement.
\begin{theorem}
Let $n=3$, and assume (\ref{ea6}) (\ref{ea7}). Consider the problem (\ref{1.1}) with $p >2$. Then there exists a global solution $u$ for any initial data which are sufficiently small, decaying and regular.
\end{theorem}
Let us describe the strategy of the proof. We basically follow the approach that appeared in \cite{W15} to give the proof. Firstly, by local existence (Lemma \ref{4.8}), we obtain for any $T > 0$ there exists $\ep_{0}$ such that if the norm of initial data less than $\ep_{0}$, the solution could exist up to time $T$. Thus from $T$ time, we can convert the damping term $\frac{\mu}{(1+t)^{\beta}}$ to some $\frac{\tilde{\mu(t)}}{(1+t)^{\tilde{\beta}}}$ with $|\tilde{\mu}|$ small enough and $\tilde{\beta}>1$ and then we mainly exploit local energy estimate with variable coefficient (Lemma \ref{7.4}) appeared in \cite{W15} Lemma 3.5 to get global solution.
\begin{remark}
The similar argument works for nonlinearity $c_{1}(u)u^{2}_{t} + c_{2}(u)|\nabla u|^{p}$, where $c_{1}, c_{2}$ are given smooth functions.
\end{remark}
\begin{remark}
When $n \geq 4$ and $g_{2}=0$, the argument in \cite{W15} can be adapted to prove the global existence for $p > p_{c}$ with small radial data.
\end{remark}
\begin{remark}
When $n=2$ and $g_{1}=g_{2}=0$, we can use the similarly argument to show (\ref{1.1}) admits global solution for $p > p_{c}$.
\end{remark}
\begin{remark}
When $n\geq 2$ and $p=p_{c}$, we can obtain the lower bound of existence time of (\ref{1.1}) $T_{\ep} \geq e^{c\ep^{1-p}}$ in a similar way.
\end{remark}
We close this section by listing some notations. The vector fields to be used will be labeled as
\[Y=(Y_1,\dots, Y_{n(n+1)/2}) = \{\nabla_x,\Omega\},\]
Here $\Omega$ denotes the
generators of spatial rotations:
\[\Omega_{ij} = x_i\partial_j - x_j\partial_i,\quad 1\le i<j\le n.\]
For a norm $X$ and a nonnegative integer $m$, we shall use the shorthand
\[|Y^{\le m} u| = \sum_{|\mu|\le m} |Y^\mu u|,\quad \|Y^{\le m} u\|_X
= \sum_{|\mu|\le m} \|Y^\mu u\|_X,\]
For fixed $T>0$, the space-time norm $L_T^qL_x^r$ is simply $L_t^q([0,T],L_x^r(\R^n))$. In the case of $T=\infty$, we use $L_t^qL_x^r$ to denote $L_t^q([0,\infty),L_x^r(\R^n))$. As usual, we use $\|\cdot\|_{E_m}$ to denote the energy norm of order $m\ge0$,
\begin{align}
\label{ea8}
&\|u\|_E=\|u\|_{E_0}=\|\partial u\|_{L_T^\infty L_x^2},\ \|u\|_{E_m}=\sum_{|a|\leq m}\|Y^au\|_E.
\end{align}
We will use $\|\cdot\|_{LE}$ to denote the (strong) local energy norm
\begin{align}
\label{ea9}
&\|u\|_{LE}=\|u\|_E+\|\partial u\|_{l_\infty^{-1/2}(L_t^2L_x^2)}+\|r^{-1}u\|_{l_\infty^{-1/2}(L_t^2L_x^2)},\ \|u\|_{LE_m}=\sum_{|a|\leq m}\|Y^au\|_{LE},
\end{align}
where we write
$$\|u\|_{l_{q}^s(A)}=\|(\phi_j(x)u(t,x))\|_{l^{s}_q(A)},$$
for a partition of unity subordinate to the dyadic (spatial) annuli, $\sum_{j\ge 0}\phi_j^2(x)=1$.

\section{Preliminary}
In this section, we give some energy type estimates for linear damped wave equation
\begin{equation}
\label{ea1}
u_{tt}-\Delta_\gm{u} + \frac{\mu \pt{u}}{(1 + t)^\be}  = F.
\end{equation}
Moreover, we list the local energy estimate Lemma \ref{7.4} which we shall use later.
\begin{lemma}
\label{74}
Assume (\ref{ea6}) (\ref{ea7}). Let $u \in C^{2}([0,T]\times \R^{n})$ vanishing for large $x$ and satisfy (\ref{ea1}) then we have
\begin{equation}
\label{1d3}
\|\partial u(t)\|_{L_x^2}\leq C_0\left(\|\partial u(0)\|_{L_x^2}+\int_0^T\|F\|_{L_x^2}dt\right),
\end{equation}
where $C_{0}$ is constant depends on $\mu$ and $\beta$.
\end{lemma}
\begin{prf}
This is standard energy type estimate (See, e.g., Sogge \cite{sogge} Proposition 2.1). For the readers' convenience, we give the sketch of proof.\par
Let $E^{2}(t) = \frac{1}{2}\int_{\R^{n}}\big((\pa_{t}u)^{2}+g^{jk}\pa_{x_{j}}u\pa_{x_{k}}u\big)\sqrt{|\gm|}dx$, then
\begin{align*}
(E^{2})'&= \frac{d}{dt}\frac{1}{2}\int_{\R^{n}}\big((\pa_{t}u)^{2}+g^{jk}\pa_{x_{j}}u\pa_{x_{k}}u\big)\sqrt{|\gm|}dx\\
&=\int_{\R^{n}}(u_{tt}-\Delta_{\gm}u)u_{t}\sqrt{|\gm|}dx\\
&=\int_{\R^{n}}Fu_{t}\sqrt{|\gm|}dx - \int_{\R^{n}}\frac{\mu}{(1+t)^{\beta}}u_{t}^{2}\sqrt{|\gm|}dx\\
&\les \|F\|_{L^{2}}E + \frac{|\mu|}{(1+t)^{\beta}}E^{2}.
\end{align*}
Thus $$E' \les \|F\|_{L^{2}}+ \frac{|\mu|}{(1+t)^{\beta}}E, $$
by Gronwall's inequality
$$E(t) \les (E(0) + \|F\|_{L^{2}})e^{\frac{|\mu|}{\beta-1}},~0\leq t \leq T.$$
By the assumption (\ref{ea6}) (\ref{ea7}), there exists a constant $C$ such that
$$\frac{1}{C}\|\pa u\|_{L^{2}} \leq E(t) \leq C\|\pa u \|_{L^{2}},$$
hence (\ref{1d3}) follows.\par
\end{prf}
\begin{corollary}(Higher order energy estimate)
\label{2d2}
Under the same assumption of Lemma \ref{74}. If $\beta > 1$, we have
\begin{equation*}
\|\partial Y^{\leq2}u(t)\|_{L_x^2}\leq C_1e^{C_1T}\left(\|\partial Y^{\leq2}u(0)\|_{L_x^2}+\int_0^T\|Y^{\leq2}F\|_{L_x^2}dt\right) ,
\end{equation*}
where $C_{1}$ is constant depends on $\mu$ and $\beta$.
\end{corollary}
\begin{prf}
Applying the vector field $Y^{\leq2}$ to both sides of the equation (\ref{ea1}), and note that
$$[\Delta_{\gm}, \nabla^{\mu}\Omega^{\nu}]u = \sum_{\substack{|\tilde{\mu}|+|\tilde{\nu}|\le
    |\mu|+|\nu|\\|\tilde{\nu}|\le |\nu|}} \tilde{b}^\alpha_{\tilde{\mu}\tilde{\nu}} \nabla^{\tilde{\mu}} \Omega^{\tilde{\nu}} \nabla_\alpha u.$$
By assumption (\ref{ea6}) (\ref{ea7}) we have $\|\tilde{b}\|_{L^{\infty}} < \infty$. Then we have
\begin{equation}
\label{2dd}
(\pa_{tt} -\Delta_{\gm})Y^{\leq2}u + \frac{\mu \pt Y^{\leq2}u}{(1 + t)^\be}  = \widetilde{F}:=Y^{\leq2}F + \sum_{|\mu|+|\nu|\leq2}[\Delta_{\gm}, \nabla^{\mu}\Omega^{\nu}]u
\end{equation}
Applying Lemma \ref{74} to (\ref{2dd}) we obtain
\begin{align*}
&\|\partial Y^{\leq2}u(t)\|_{L_x^2}\\
\leq &C_0\big(\|\partial Y^{\leq2}u(0)\|_{L_x^2}+\int_0^T\|\widetilde{F}\|_{L_x^2}dt\big) \\
\les &\big(\|\partial Y^{\leq2}u(0)\|_{L_x^2}+\int_0^T\|Y^{\leq2}F\|_{L_x^2}dt +\sum_{|\mu|+|\nu|\leq2}\int_0^T\|[\Delta_{\gm}, \nabla^{\mu}\Omega^{\nu}]u\|_{L_x^2}dt\big)\\
\les&\big(\|\partial Y^{\leq2}u(0)\|_{L_x^2}+\int_0^T\|Y^{\leq2}F\|_{L_x^2}dt\big)+ \int_{0}^{T}\|\tilde{b}\|_{L_{x}^{\infty}}\|\pa Y^{\leq 2}u\|_{L^{2}_{x}}dt
\end{align*}
Then by Gronwall's inequality,
\begin{align*}
\|\partial Y^{\leq2}u(t)\|_{L_x^2}\leq C_1e^{C_1T}\left(\|\partial Y^{\leq2}u(0)\|_{L_x^2}+\int_0^T\|Y^{\leq2}F\|_{L_x^2}dt\right).
\end{align*}
\end{prf}

\begin{lemma}
\label{7.4}
Let $n\geq3$ and consider the problem $u_{tt} - \Delta_{\gm}=F$ on the manifold $(\R^{n},\gm)$ satisfying (\ref{ea6}) (\ref{ea7}) with $\rho=min(\rho_1, \rho_2-1)>0$. Then for any positive $\theta>0$, we have the following higher order local energy estimates:
\begin{eqnarray}
\|u\|_{LE_k}\les \sum_{|a|\les k}\left(\|\partial Y^au(0)\|_{L_x^2}+\|Y^aF\|_{L_t^1L_x^2+l_2^{\frac12+\theta}L_t^2L_x^2}\right).
\end{eqnarray}
\end{lemma}
\begin{prf}
 See Wang \cite{W15} Lemma 3.5.
\end{prf}

\section{Local Existence}
\begin{lemma}
 \label{4.8}
Consider the problem (\ref{1.1}) on the manifold $(\R^{3}, \gm)$ satisfying (\ref{ea6}) (\ref{ea7}). Then there exists small positive constant $\ep_{0}$, such that for any initial data $(f,g) \in H^{1} \times L^{2}$ satisfying $$\|\nabla Y^{\leq2}f\|_{L_x^2}+\|Y^{\leq2}g\|_{L_x^2}=\ep^p \leq \ep_{0},$$
there exist $T_{\ep}>0$ with $\lim_{\ep \rightarrow 0^{+}}T_{\ep} = \infty$. and a unique solution $u$ with
$$\|\partial Y^{\leq2}u(t)\|_{L_t^{\infty}([0,T])L_x^2}\leq2 C_1\ep.$$
\end{lemma}
\begin{prf}
We define a closed subset of Banach space
$$E_T=\{u;\ \|u\|_{E_T}:=\sup_{0\leq t\leq T}\|\partial Y^{\leq2}u\|_{L_x^2}\leq2 C_1\ep\},$$
and a map $\Ga:v\to u$, satisfying
\beeq
\left\{\begin{array}{l}u_{tt} - \Delta_{\gm}u+\frac{\mu\partial_{t} u}{(1+t)^{\be}}=|\pt v|^p,\quad
\\
u(0,x) = f,
 \pt u(0,x) = g.\end{array}
 \right.
 \label{3.1}
\eneq
Then we are reduced to show $\Ga$ is a contraction map from $E_T$ to $E_T$.

Firstly, we will prove $\Ga :E_T\to E_T$.\\
For any $v\in E_T$, by Corollary \ref{2d2}, we have
\begin{align*}
\|\Ga v\|_{E_T}&=\sup_{0\leq t\leq T}\|\partial Y^{\leq2}u\|_{L_x^2}\\
&\leq C_1e^{C_1T}\left(\|\partial Y^{\leq2}u(0)\|_{L_x^2}+\int_0^T\|Y^{\leq2}F\|_{L_x^2}dt\right) .
\end{align*}
Since
\begin{eqnarray*}
Y^{\leq2}F(v)=Y^{\leq2}|\pt v|^p\sim|\pt v|^{p-1}Y^{\leq2}\pt v+|\pt v|^{p-2}(Y^{\leq1}\partial v)^2.
\end{eqnarray*}
For the first term, by H$\ddot{\mathrm{o}}$lder's inequality and Sobolev imbedding $H^{2}(\R^{3}) \hookrightarrow L^{\infty}_{x}(\R^{3})$ we have
\begin{eqnarray*}
\||\pt v|^{p-1}Y^{\leq2}\pt v\|_{L_x^2}\leq\|\pt v\|_{L^{\infty}}^{p-1}\|\pa Y^{\leq2}v\|_{L_x^2} \les \|\pa Y^{\leq2}v\|^{p}_{L_x^2}.
\end{eqnarray*}
Similarly, for the second term, by Sobolev imbedding $H^{1}(\R^{3}) \hookrightarrow L^{4}(\R^{3})$, we have
\begin{eqnarray*}
\||\pt v|^{p-2}(Y^{\leq1}\partial v)^2\|_{L_x^2}\leq\|\pt v\|^{p-2}_{L_x^{\infty}}\|Y^{\leq1}\partial v\|_{L_x^4}^2 \les \|\pt v\|^{p-2}_{L_x^{\infty}}\|\pa Y^{\leq 2}v\|^{2}_{L^{2}_{x}} \les \|\pa Y^{\leq2}v\|^{p}_{L_x^2}.
\end{eqnarray*}
Hence
\begin{align*}
\|\Ga v\|_{E_T} &\leq C_1e^{C_1T}(\ep^{p} + C_{2}(2C_{1}\ep)^{p})\\
&\leq 2C_{1}\ep
\end{align*}
if we take
$$T \leq \frac{1}{C_{1}}\ln(\frac{1}{2^{p}C_{1}^{p}C_{2}\ep^{p-1}}).$$
Secondly, we will prove $\Ga :E_T\to E_T$ is contraction.
For any $u,\,v\in E_{T}$ with $(\Ga u, u)$, $(\Ga v, v)$ satisfy (\ref{3.1}), then $\Ga (u- v)$ satisfies that
\beeq
\left\{\begin{array}{l}(\pa_{tt} -\Delta_{\gm})\Ga (u- v)+\frac{\mu\pt \Ga (u- v)}{(1+t)^{\be}}=|\pt u|^p-|\pt v|^p,\quad
\\
\Ga (u-v)(0,x) = 0,
 \pt \Ga (u- v)(0,x) = 0.\end{array}
 \right.
 \label{3.2}
\eneq
By applying corollary \ref{2d2}, we have
\begin{eqnarray*}
\|\partial Y^{\leq2}\Ga (u-v)\|_{L_x^2}
\leq C_1e^{C_1T}\int_0^T\|Y^{\leq2}(|\pt u|^p-|\pt v|^p)\|_{L_x^2}dt.
\end{eqnarray*}
Since
\begin{eqnarray*}
Y^{\leq2}(|\pt u|^p-|\pt v|^p)\sim (|\pt u|^{p-3}+|\pt v|^{p-3})(Y(\pt u+\pt v))^2(\pt u-\pt v)\\
+(|\pt u|^{p-2}+|\pt v|^{p-2})Y^{\leq2}(\pt u+\pt v)Y^{\leq2}(\pt u-\pt v),
\end{eqnarray*}
by  H$\ddot{\mathrm{o}}$lder's inequality and Sobolev inequality, we have
\begin{align*}
&\||\pt u|^{p-3}(Y(\pt u+\pt v))^2(\pt u-\pt v)\|_{L_x^2}\\
\les &\|\pt u\|_{L_x^{\infty}}^{p-3}\|Y(\pt u+\pt v))\|_{L_x^4}^2\|\pt u-\pt v\|_{L_x^{\infty}}\\
\les &\|\partial Y^{\leq2}u\|_{L_x^2}^{p-3}\|Y^{\leq2}\partial u+Y^{\leq2}\partial  v))\|_{L_x^2}^2\|Y^{\leq2}\partial(u-v)\|_{L_x^2}.
\end{align*}
Similarly,
\begin{eqnarray*}
&&\|(|\pt u|^{p-2}+|\pt v|^{p-2})Y^{\leq2}(\pt u+\pt v)Y^{\leq2}(\pt u-\pt v)\|_{L_x^2}\\
&\les&(\|\pt u\|_{L_x^{\infty}}^{p-2}+\|\pt v\|_{L_x^{\infty}}^{p-2})\|Y^{\leq2}(\pt u+\pt v)\|_{L_x^4}\|Y^{\leq2}(\pt u-\pt v)\|_{L_x^4}\\
&\les&(\|\partial Y^{\leq2} u\|_{L_x^2}^{p-2}+\|\partial Y^{\leq2} v\|_{L_x^2}^{p-2})\|Y^{\leq2}\partial u+Y^{\leq2}\partial  v\|_{L_x^2}\|Y^{\leq2}\partial (u-v)\|_{L_x^2}.
\end{eqnarray*}
Hence
\begin{align*}
\|\Ga (u-v)\|_{E_{T}}&=\|\partial Y^{\leq2}(u-v)\|_{L^{\infty}_{t}L_x^2}\\
& \les C_1e^{C_1T}TC_{3}(2C_1\ep)^{p-1}\|u-v\|_{E_T}\\
& \leq \frac12\|u-v\|_{E_T},
\end{align*}
if we take
$$T\leq \frac{1}{C_{1}}\ln(\frac{1}{2^{p}C_{1}^{p}C_{3}\ep^{p-1}}).$$
Therefore, $\Ga$ is a contraction map from $E_T$ to $E_T$, if we take
 $$T=\min\{\frac{1}{C_{1}}\ln(\frac{1}{2^{p}C_{1}^{p}C_{3}\ep^{p-1}}), \frac{1}{C_{1}}\ln(\frac{1}{2^{p}C_{1}^{p}C_{2}\ep^{p-1}})\}.$$
Thus there is a unique solution in $E_{T}$.
\end{prf}

\section{Global Existence of equation (\ref{1.1})}
\begin{theorem}
\label{4.1}
Consider the problem (\ref{1.1}) on the manifold $(\R^{3},\gm)$ satisfying (\ref{ea6}) (\ref{ea7}) with $\rho=min(\rho_1, \rho_2-1)>0$. Then there exists small positive constants $\ep_0$ such that for any initial data satisfying
\beeq
\label{4.3}
 \sum_{|a|\leq2}\|\partial Y^au(0)\|_{L_x^2(\R^3)}=\ep^{p}
\leq \ep_0,\
\|u(0)\|_{L^2(\R^3)}<\infty\ .
\eneq
There is a global solution $u$ with $\|u\|_{LE_2}\les \ep$.
\end{theorem}
\begin{prf}
We rewrite the equation (\ref{1.1}) as
\begin{equation}
\label{704}
u_{tt} - \Delta_\gm{u} + \frac{\mu}{(1 + t)^{\gamma}}\frac{\pt u}{(1+t)^{\be-\gamma}}  = |\pt u|^p,
\end{equation}
where $0<\gamma < 1$ such that $\tilde{\be} = \be-\gamma>1$. Let $\tilde{\mu}(t)=\frac{\mu}{(1+t)^{\gamma}}$, then
(\ref{704}) becomes
\beeq
\label{705}
u_{tt} - \Delta_\gm{u} +  \frac{\tilde{\mu}(t)u_{t}}{(1+t)^{\tilde{\beta}}}= |\pt u|^p,
\eneq
From the local existence Lemma \ref{4.8}, we know for any $T>0$, there exits $\ep_0>0$, if initial data
satisfy (\ref{4.3}) then (\ref{705}) has a unique solution $u$ in $[0,T]$ with
\begin{equation*}
\|\partial Y^{\leq2}u\|_{L_t^{\infty}([0,T])L_x^2}\leq2C_1\ep.
\end{equation*}
Thus $|\tilde{\mu}(t)|\leq\frac{|\mu|}{(1+T)^{\gamma}}$ can be sufficiently small if we take $T$ sufficiently large. Then from the $T$ time, we are reduced to consider the equation
\beeq
\begin{cases}
u_{tt} - \Delta_\gm{u}=G(u):=|\pt u|^p-\frac{\tilde{\mu}(t)\pt u(t)}{(1+t)^{\tilde{\be}}}, t \geq T\\
u(T,x) = f, \pt u(T,x) = g.
\end{cases}
\eneq
with
\begin{align}
\label{706}
&\|Y^{\leq 2}\nabla f\|_{L_x^2} + \|Y^{\leq 2}g\|_{L_x^2}\leq 2C_{1}\ep,
\end{align}
and
$$\|u(T)\|_{L^{2}} = \|\int_{0}^{T}\pa_{t}u(s)ds + u(0)\|_{L^{2}} \leq T\|\pa u\|_{L^{\infty}_{t}L^{2}_{x}}+\|u(0)\|_{L^{2}} < \infty.$$
Let $G(u) = G_{1}(u) + G_{2}(u) = |u_{t}|^{p} - \frac{\tilde{\mu}(t)\pt u(t)}{(1+t)^{\tilde{\be}}}$ and
from now on, we denote $\|u\|_{L_t^pL_x^q}$ as $\|u\|_{L_t^p([T,{\infty}])L_x^q}$.
Set $u_0\equiv0$ and define $u_{k+1}$ to be the solution to the equation
$$(\pa_{tt} - \Delta_{\gm})u_{k+1}=G(u_k),u_{{k+1}}(T,x)=f(x),\pt u_{{k+1}}(T,x)=g(x).$$

$Boundedness$: By the smallness condition (\ref{706}) on the data, it follows from Lemma \ref{7.4} that there is a universal constant $C_4$ so that
\begin{align}
\label{ea4}
&\|u_1\|_{LE_2}\leq C_4\ep, \ \|u_{k+1}\|_{LE_2}\leq C_4\ep+C_4\sum_{|a|\leq 2}\|Y^aG(u_k)\|_{L_t^1L_x^2}.
\end{align}
We shall argue inductively to prove that
\begin{align}
\label{ea5}
&\|u_{k+1}\|_{LE_2}\leq 2C_4\ep.
\end{align}
By the above, it suffices to show
\begin{align}
\label{ea60}
&\sum_{|a|\leq 2}\|Y^aG(u)\|_{L_t^1L_x^2}\leq \ep,
\end{align}
for any $u$ with $\|u\|_{LE_2}\leq 2C_4\ep\leq 1.$\\
For the $G_{1} = |u_{t}|^{p}$ part, by (4.6) of \cite{W15}
\begin{equation*}
\|Y^{\leq 2}|\pa_{t} u|^p\|_{L_t^1L_x^2}\les \|u\|_{LE_2}^p.
\end{equation*}
For the $G_{2}(u)$ part, since
\begin{align*}
\|\frac{\tilde{\mu}(t)}{(1+t)^{\tilde{\be}}}Y^{\leq 2}\pt u\|_{L_t^1L_x^2}
&\les \|\frac{\tilde{\mu}(t)}{(1+t)^{\tilde{\be}}}\|_{L_t^1}\|Y^{\leq 2}\pt u\|_{L_t^{\infty}L_x^2}\\
&\les\frac{\|\tilde{\mu}(t)\|_{L_t^{\infty}}}{\tilde{\be}-1}\|u\|_{LE_2}
\end{align*}
In conclusion, we see that there exists a constant $C_{5}$ such that
$$\|Y^{\leq 2}G(u)\|_{L_t^1L_x^2}\leq C_{5}(2C_{4}\ep)^{p} + 2C_{4}\ep C_{5}\frac{\|\tilde{\mu}(t)\|_{L_t^{\infty}}}{\tilde{\be}-1} \leq \ep$$
for $\ep \leq \ep_{0}$ with
$$C_{5}(2C_{4})^{p}\ep_{0}^{p-1} + 2C_{4} C_{5}\frac{\|\tilde{\mu}(t)\|_{L_t^{\infty}}}{\tilde{\be}-1} \leq 1.$$
$Convergence\ of\ the\ sequence\ {u_k}$:
We see that
$$\|G(u)-G(v)\|_{L_t^1L_x^2}\leq \|G_1(u)-G_1(v)\|_{L_t^1L_x^2} + \|G_2(u)-G_2(v)\|_{L_t^1L_x^2}.$$
For the first part, by the 1st display on page 7442 of \cite{W15} we have
$$\|G_1(u)-G_1(v)\|_{L_t^1L_x^2} \les (\|u\|_{LE_{2}}+\|v\|_{LE_{2}})^{p-1}\|u-v\|_{LE}.$$
For the second part, since
\begin{align*}
\|G_2(u)-G_2(v)\|_{L_t^1L_x^2}&= \|\frac{\tilde{\mu}(t)}{(1+t)^{\tilde{\be}}}(\pt u-\pt v)\|_{L_t^1L_x^2}\\
&\les \|\frac{\tilde{\mu}(t)}{(1+t)^{\tilde{\be}}}\|_{L_t^1L_x^{\infty}}\|\partial(u-v)\|_{L_t^{\infty}L_x^2}\\
&\les \frac{\|\tilde{\mu}(t)\|_{L_t^{\infty}}}{\tilde{\be}-1}\|u-v\|_{LE}.
\end{align*}
Hence there exists a constant $C_{6}$ such that
$$\|G(u)-G(v)\|_{L_t^1L_x^2} \leq \big(C_{6}(4C_{4}\ep)^{p-1}+ C_{6}\frac{\|\tilde{\mu}(t)\|_{L_t^{\infty}}}{\tilde{\be}-1}\big)\|u-v\|_{LE}\leq \frac{1}{2}\|u-v\|_{LE}$$
for $\ep \leq \ep_{0}$ and $\ep_{0}$ satisfies
$$\big(C_{6}(4C_{4}\ep_{0})^{p-1}+ C_{6}\frac{\|\tilde{\mu}(t)\|_{L_t^{\infty}}}{\tilde{\be}-1}\big)\leq \frac{1}{2}.$$
Together with the uniform boundedness (\ref{ea5}), we find an unique global solution $u \in L^{\infty}_{x}([T,\infty];H^{3})\cap Lip_{t}([T,\infty];H^{2})$ with $\|u\|_{LE_{2}}\leq 2C_{4}\ep$.
\end{prf}
\section*{acknowledgements}
The authors would like to thank their advisor Professor Chengbo Wang, for introducing the problem and kind support. Without his consistent and illuminating instruction, this paper could not have reached its present stage.
\bibliographystyle{plain1}

\end{document}